\documentclass[11pt]{article} \topmargin=-0.5cm
\usepackage{amssymb}

   \csname @addtoreset\endcsname{equation}{section}
    \csname @addtoreset\endcsname{figure}{section}
    \csname @addtoreset\endcsname{table}{section}

\newcounter{thanksnum}
\def\thanksnumber#1
{\setcounter{thanksnum}{\value{footnote}}\setcounter{footnote}{#1}%
                     \addtocounter{footnote}{-1}\footnotemark
                     \setcounter{footnote}{\value{thanksnum}}}
\def\newtheoremz#1{\@ifnextchar[{\@othmz{#1}}{\@nthmz{#1}}}

\def\@nthmz#1#2{%
\@ifnextchar[{\@xnthmz{#1}{#2}}{\@ynthmz{#1}{#2}}}

\def\@xnthmz#1#2[#3]{\expandafter\@ifdefinable\csname #1\endcsname
{\@definecounter{#1}\@addtoreset{#1}{#3}%
\expandafter\xdef\csname the#1\endcsname{\expandafter\noexpand
  \csname the#3\endcsname \@thmcountersepz \@thmcounterz{#1}}%
\global\@namedef{#1}{\@thmz{#1}{#2}}\global\@namedef{end#1}{\@endtheoremz}}}

\def\@ynthmz#1#2{\expandafter\@ifdefinable\csname #1\endcsname
{\@definecounter{#1}%
\expandafter\xdef\csname the#1\endcsname{\@thmcounterz{#1}}%
\global\@namedef{#1}{\@thm{#1}{#2}}\global\@namedef{end#1}{\@endtheoremz}}}

\def\@othmz#1[#2]#3{\expandafter\@ifdefinable\csname #1\endcsname
  {\global\@namedef{the#1}{\@nameuse{the#2}}%
\global\@namedef{#1}{\@thmz{#2}{#3}}%
\global\@namedef{end#1}{\@endtheoremz}}}

\def\@thmz#1#2{\refstepcounter
    {#1}\@ifnextchar[{\@ythmz{#1}{#2}}{\@xthmz{#1}{#2}}}

\def\@xthmz#1#2{\@begintheoremz{#2}{\csname the#1\endcsname}\ignorespaces}
\def\@ythmz#1#2[#3]{\@opargbegintheoremz{#2}{\csname
       the#1\endcsname}{#3}\ignorespaces}

\def\@thmcounterz#1{\noexpand\arabic{#1}}
\def\@thmcountersepz{.}
\def\@begintheoremz#1#2{ \trivlist \item[\hskip \labelsep{\bf #1\ #2}]}
\def\@opargbegintheoremz#1#2#3{ \trivlist
      \item[\hskip \labelsep{\bf #1\ #2\ (#3)}]}
\def\@endtheoremz{\endtrivlist}

\newtheorem{theorem}{Theorem}[section]
\newtheorem{lemma}{Lemma}[section]

\newtheorem{proposition}{Proposition}[section]
\newtheorem{corollary}{Corollary}[section]
\newtheorem{condition}{Condition}[section]
\newtheorem{definition}{Definition}[section]
\newtheorem{remark}{Remark}[section]

\newtheorem{example}{Example}[section]

\def\e{\varepsilon}

\def\defi{\stackrel{{\scriptscriptstyle \Delta}}{=}}

\def\d{\delta}
\def\o{\omega}
\def\O{\Omega}
\def\q{q}

\def\F{{\cal F}}
\def\w{\widehat}
\def\Ind{{\mathbb{I}}}

\def\mes{{\rm mes\,}}

\def\esssup{\mathop{\rm ess\, sup}}

\def\R{{\bf R}}
\def\E{{\bf E}}
\def\P{{\bf P}}

\def\L{L}

\def\b{\beta}
\def\s{\delta}
\def\g{\gamma}
\def\C{{\bf C}}

\def\ww{\widetilde}
\def\X{{\cal X}}
\def\t{\theta}
\def\oo{\bar}
\def\s{\sigma}
\def\p{\partial}
\def\G{\Gamma}

\def\A{{\cal A}}
\def\M{{\cal M}}

\def\L{{\cal L}}

\def\TT{{\cal T}}
\newcommand{\be}{\begin{equation}}
\newcommand{\ee}{\end{equation}}
\newcommand{\bd}{\begin{displaymath}}
\newcommand{\ed}{\end{displaymath}}
\newcommand{\ba}{\begin{array}{ll}}
\newcommand{\ea}{\end{array}}
\newcommand{\baa}{\begin{eqnarray}}
\newcommand{\eaa}{\end{eqnarray}}
\newcommand{\baaa}{\begin{eqnarray*}}
\newcommand{\eaaa}{\end{eqnarray*}}   \font\sm=cmr10
\def\k{\kappa}

\def\ww{\tilde}

\def\QQ{{ Q}}
\def\Q{{\cal Q}}
\def\CC{{\cal C}}
\title{On almost surely
periodic and almost periodic solutions of backward SPDEs }
\author{
Nikolai Dokuchaev\\
 {\sm Department of Mathematics \& Statistics, Curtin
University,}\\ {\sm  GPO Box U1987, Perth, 6845 Western Australia} \\ {\sm Email N.Dokuchaev@curtin.edu.au}}
\begin{document}
\maketitle
\begin{abstract}
We study linear  backward stochastic partial differential equations
of parabolic type with special boundary conditions in time. The
standard Cauchy condition at the terminal time is replaced by a
condition that holds almost surely and mixes the random values of
the solution at different times, including the terminal time,
initial time and continuously distributed times. Uniqueness,
solvability and regularity results for the solutions are obtained.
In particular, conditions of existence of periodic in time and
"almost periodic" solutions are obtained for backward SPDEs.
\end{abstract}
%{\it Abbreviated tittle: }
\section{Introduction}
Partial differential equations and  stochastic partial
differential equations (SPDEs) have fundamental significance for
natural sciences, and various boundary value problems for them
were widely studied.
 Usually,  well-posedness of a  boundary value depends on
the choice of the boundary value conditions. For the deterministic
parabolic equations, well-posedness
 requires  the correct choice
 of the initial condition. For example, consider  the heat equation $u'_t=u''_{xx}$,
$t\in[0,T]$. For this equation, a boundary value problem with the
Cauchy condition at initial time $t=0$ is well-posed, and a
boundary value problem with the Cauchy condition at terminal time
$t=T$ is ill-posed.
 It is known also that the problems for deterministic  parabolic
equation are well-posed for  periodic type condition
$u(x,0)=u(x,T)$; see, e.g.,  Dokuchaev (1994, 1995), Fife (1964),
Hess (1991), Lieberman  (1999), Nakao (1984),  Shelukhin (1993),
Vejvoda (1982).  Less is known for parabolic equation with more
general non-local in time conditions and for SPDEs.
\par
 Boundary value problems for SPDEs are well studied  in the existing literature
 for the case
 of  forward parabolic Ito equations with the  Cauchy condition at
initial time (see, e.g., Al\'os et al (1999), Bally {\it et al}
(1994),  Da Prato and Tubaro (1996), Gy\"ongy (1998), Krylov (1999),
Maslowski (1995), Pardoux (1993),
 Rozovskii (1990), Walsh (1986), Zhou (1992),
and the bibliography there).  Many results have been also obtained
for the backward
 parabolic Ito equations with  Cauchy
 condition   at terminal time, as well as for  pairs of   forward and
backward
 equations with separate  Cauchy conditions at initial time and
 the terminal time respectively; see, e.g., Yong and Zhou
 (1999), and the author's papers (1992), (2005), (2011), (2012).  Note that a
 backward SPDE cannot be transformed into a forward
  equation by a simple time
 change, unlike as for the case of deterministic equations. Usually, a
  backward SPDE is solvable in the sense that there exists a
  diffusion term being considered as a part of the solution that
  helps to ensure that the solution is adapted to the driving Brownian
motions.
\par
There are also results for SPDEs with boundary conditions that mix
the solution at different times that may  include initial time and
terminal time. This category includes stationary type solutions for
 forward SPDEs (see, e.g., Dorogovtsev and  Ortega (1988), Caraballo {\em et al } (2004),
 Chojnowska-Michalik (19987), Chojnowska-Michalik and Goldys
 (1995), Duan {\em et al} (2003), Mattingly (1999),
Mohammed  {\em et al} (2008), Sinai (1996), and the references
here). There are also results for different types of the periodicity
of the solutions of SPDEs; see, e.g., Chojnowska-Michalik (1990),
Tudor (1992), Da Prato and  Tudor (1995), Arnold and Tudor (1998),
Kl\"unger  (2001),  Bezandry and Diagana (2007),  Mellah and de Fitte (2007), Feng and Zhao
(2012), Bedouhene {\em et al}  (2012),  Crewe  (2013).
 As was mentioned in Feng and Zhao (2012), it is difficult to expect that,
 in general, a SPDE has a periodic
 in time solution $u(\cdot,t)|_{t\in[0,T]}$ in a usual sense of exact equality
 $u(\cdot,t)=u(\cdot,t+T)$ that holds almost surely
 given that $u(\cdot,t)$ is adapted to some Brownian motion.
 However, there are important
 examples of stochastic processes with this property. In particular, this property  holds at  $t=0$  for a
 Brownian bridge. Using this, the existence of almost surely periodic
solutions  was established in
Rodkina (1992) for ordinary stochastic equations  with the
driving Brownian motion replaced by a Brownian bridge.

\par
In a more typical setting with driving Brownian motion, the
periodicity of the solutions of stochastic equations has to be
interpreted differently.
 This  periodicity was usually considered
 in the sense of the distributions. In Feng and Zhao (2012),
 the periodicity was established  in a  stronger sense as a random
 periodic solution (see Definition 1.1 from Feng and Zhao (2012));
 this definition does not assume the equality  $u(\cdot,t)=u(\cdot,T)$.
In Feng and Zhao (2012),  semi-linear parabolic Ito equations with a
self-adjoint main operator were considered. There are also results
for  almost periodic in mean-square sense solutions; see, e.g.,
Tudor (1992), Da Prato and  Tudor (1995), Arnold and Tudor (1998),
Bezandry and Diagana (2007),  Bedouhene {\em et al}  (2012), Crewe
(2013).

The present paper addresses these and related problems again for a single period setting.
We found examples of SPDEs  where almost periodicity conditions and exact periodicity conditions hold almost surely a
well as more general non-local boundary value conditions. It appears
that this is possible with the replacement of forward SPDEs for
backward SPDEs.

We consider linear  Dirichlet  condition at the boundary of the
state domain; the equations are of  a parabolic type  and are not
necessary self-adjoint. The standard boundary value Cauchy condition
at the one fixed time is replaces by a condition that mixes in one
equation the values of the solution at different times over given
time interval, including the terminal time and continuously
distributed times. This is a novel setting comparing with the
periodic conditions for the distributions, or with  conditions from
Kl\"unger (2001) and Feng and Zhao (2012), or with conditions for
the expectations from Dokuchaev (2008), or mean-square almost
periodicity from Tudor (1992). These conditions include, for
instance, conditions $\kappa u(\cdot,0)=u(\cdot,T)$ a.e. with
$\kappa\in\R$  (Theorems \ref{Th1}-\ref{Th5}). We present sufficient
conditions for existence and regularity of solutions  in
$L_2$-setting  (Corollary \ref{corrAP} and Corollary \ref{corrAP2}).
Corollary \ref{corrAP2} establishes  existence of modified almost
periodic solution. It can be noted that we consider "almost
periodicity" on a single time period; this setting is easier than a
multiperiod setting on infinite interval since it does not require
to formalize a periodic extension of the coefficients. However, with
respect to a single period, our "almost periodicity" property from
Corollary \ref{corrAP2}) is stronger than the mean square "almost
periodicity". Finally, we present sufficient conditions for
existence and regularity of almost surely exact periodicity (Theorem
\ref{Th5}). The proofs is based on compactness and Fredholm theory
in $L_2$-spaces.
 The periodic solution
obtained here can be considered as a generalization of a
classical
 Brownian bridge for the case of an infinite dimensional state space.

It can be noted  the almost surely periodicity was achieved
for the purely backward SPDEs; the previous result of this kind was obtained for forward-backward
SPDEs in Feng and Zhao (2012). Usually, backward SPDEs are usually associated with a
Cauchy condition at the terminal time.   With regards to the general theory of SPDEs, our results open a way
to extend applications of backward SPDEs on the problems with
periodic and mixed in time conditions.

Related problems  were considered in
 Dokuchaev (2012b,c)  for a less general backward SPDE with $\oo\b_i=0$, in the notations of the present paper.   In Dokuchaev (2012b), the approach was based 
 on the contraction mapping theorem in a $L_\infty$-space;
 this approach is not applicable for  the more general SPDEs considered in the present paper. 
 In Dokuchaev (2012c), related forward and backward SPDEs with $\oo\b_i=0$ were studied in an unified framework.

 \section{The problem setting and definitions}
 We
are given a standard  complete probability space $(\O,\F,\P)$ and
a right-continuous filtration $\F_t$ of complete $\s$-algebras of
events, $t\ge 0$. We are given also a $N$-dimensional Wiener
process $w(t)$ with independent components;  it is a Wiener
process with respect to $\F_t$.
\par
Assume that we are given an open domain $D\subset\R^n$ such that
either $D=\R^n$ or $D$ is bounded  with $C^2$-smooth boundary $\p
D$. Let $T>0$ be given, and let $Q\defi D\times [0,T]$. \par
 We will study the following boundary value
problem in $Q$
\begin{eqnarray} %4.1
\label{parab1} &&d_tu+(\A u+ \varphi)\,dt +\sum_{i=1}^N
B_i\chi_idt=\sum_{i=1}^N\chi_i(t)dw_i(t), \quad t\ge 0,
\\\label{parab10}
&& u(x,t,\o)\,|_{x\in \p D}=0
\\ &&u(\cdot, T)-\G u(\cdot)=\xi.
\label{parab2}%3.1
\end{eqnarray}
Here $u=u(x,t,\o)$, $\varphi=\varphi(x,t,\o)$,
$\chi_i=\chi_i(x,t,\o)$,
 $(x,t)\in Q$,   $\o\in\O$.\par
   In (\ref{parab2}), $\G$ is a linear operator that maps functions
defined on $Q\times \O$  to functions defines on $D\times \O$. For
instance, the case where  $\G u=u(\cdot,0)$ is not excluded; this
case corresponds to the periodic type boundary condition \baa
u(\cdot,T)-u(\cdot,0)=\xi.\label{period}\eaa
\par
  In (\ref{parab1}),
  \be\label{A}\A v\defi
 \sum_{i=1}^n\frac{\p }{\p x_i}
\sum_{j=1}^n \Bigl(b_{ij}(x,t,\o)\frac{\p v}{\p x_j}(x)\Bigr)
   +\sum_{i=1}^n f_{i}(x,t,\o)\frac{\p v}{\p x_i }(x)
-\,\lambda(x,t,\o)v(x), \ee where $b_{ij}, f_i, x_i$ are the
components of $b$,$f$, and $x$ respectively, and \be\label{B}
B_iv\defi\frac{dv}{dx}\,(x)\,\beta_i(x,t,\o) +\oo
\beta_i(x,t,\o)\,v(x),\quad i=1,\ldots ,N. \ee
\par
We assume that the functions $b(x,t,\o):
\R^n\times[0,T]\times\O\to\R^{n\times n}$, $\b_j(x,t,\o):
\R^n\times[0,T]\times\O\to\R^n$, $\oo\b_i(x,t,\o):$
$\R^n\times[0,T]\times\O\to\R$, $f(x,t,\o):
\R^n\times[0,T]\times\O\to\R^n$, $\lambda(x,t,\o):
\R^n\times[0,T]\times\O\to\R$,   $\chi_i(x,t,\o): \R^n\times
[0,T]\times\O\to\R$, and $\varphi (x,t,\o): \R^n\times
[0,T]\times\O\to\R$ are progressively measurable with respect to
$\F_t$ for all $x\in\R^n$, and the function $\xi(x,\o):
\R^n\times\O\to\R$ is $\F_0$-measurable for all $x\in\R^n$. In fact,
we will also consider $\varphi$ and $\xi$ from wider classes. In
particular, we will consider generalized functions $\varphi$.
\par
If the functions $b$, $f$, $\lambda$, $\varphi$, $\G$, and $\xi$,
 are deterministic, then  $\chi_i\equiv 0$ and equation (\ref{parab1}) is
deterministic.
\subsection*{Spaces and classes of functions} %2
We denote by $\|\cdot\|_{ X}$ the norm in a linear normed space
$X$, and
 $(\cdot, \cdot )_{ X}$ denote  the scalar product in  a Hilbert space $
X$.
\par
We introduce some spaces of real valued functions.
\par
 Let $G\subset \R^k$ be an open
domain, then ${W_q^m}(G)$ denote  the Sobolev  space of functions
that belong to $L_q(G)$ together with the distributional
derivatives up to the $m$th order, $q\ge 1$.
\par
 We denote  by $|\cdot|$ the Euclidean norm in $\R^k$, and $\bar G$ denote
the closure of a region $G\subset\R^k$.
\par Let $H^0\defi L_2(D)$,
and let $H^1\defi \stackrel{\scriptscriptstyle 0}{W_2^1}(D)$ be the
closure in the ${W}_2^1(D)$-norm of the set of all smooth functions
$u:D\to\R$ such that  $u|_{\p D}\equiv 0$. Let $H^2=W^2_2(D)\cap
H^1$ be the space equipped with the norm of $W_2^2(D)$. The spaces
$H^k$ and $W_2^k(D)$ are called  Sobolev spaces, they are Hilbert
spaces, and $H^k$ is a closed subspace of $W_2^k(D)$, $k=1,2$.
\par
 Let $H^{-1}$ be the dual space to $H^{1}$, with the
norm $\| \,\cdot\,\| _{H^{-1}}$ such that if $u \in H^{0}$ then
$\| u\|_{ H^{-1}}$ is the supremum of $(u,v)_{H^0}$ over all $v
\in H^1$ such that $\| v\|_{H^1} \le 1 $. $H^{-1}$ is a Hilbert
space.
\par We shall write $(u,v)_{H^0}$ for $u\in H^{-1}$
and $v\in H^1$, meaning the obvious extension of the bilinear form
from $u\in H^{0}$ and $v\in H^1$.
\par
We denote by $\oo\ell _{k}$ the Lebesgue measure in $\R^k$, and we
denote by $ \oo{{\cal B}}_{k}$ the $\sigma$-algebra of Lebesgue
sets in $\R^k$.
\par
We denote by $\oo{{\cal P}}$  the completion (with respect to the
measure $\oo\ell_1\times\P$) of the $\s$-algebra of subsets of
$[0,T]\times\O$, generated by functions that are progressively
measurable with respect to $\F_t$.
\par
 We  introduce the spaces
 \baaa
 &&X^{k}(s,t)\defi L^{2}\bigl([ s,t ]\times\Omega,
{\oo{\cal P }},\oo\ell_{1}\times\P;  H^{k}\bigr), \quad\\ &&Z^k_t
\defi L^2\bigl(\Omega,{\cal F}_t,\P; H^k\bigr),\\
&&\CC^{k}(s,t)\defi C\left([s,t]; Z^k_T\right), \qquad k=-1,0,1,2,
\\&& \X^k_c= L^{2}\bigl([ 0,T ]\times\O,\, \oo{{\cal P}
},\oo\ell_{1}\times\P;\; C^k(\oo D)\bigr),\quad k\ge 0. \eaaa
%$\V^{k}(s,T)=L^{2}\bigl([s,T ],\oo\ell_{1}\times\P;W^{k}_2(D)\bigr)$,
%$ k=0,2,..$,
  The
spaces $X^k(s,t)$ and $Z_t^k$  are Hilbert spaces.
 \par
We introduce the spaces $$ Y^{k}(s,t)\defi
X^{k}(s,t)\!\cap \CC^{k-1}(s,t), \quad k=1,2, $$ with the norm $ \|
u\| _{Y^k(s,T)}
\defi \| u\| _{{X}^k(s,t)} +\| u\| _{\CC^{k-1}(s,t)}. $
For brevity, we shall use the notations
 $X^k\defi X^k(0,T)$, $\CC^k\defi \CC^k(0,T)$,
and  $Y^k\defi Y^k(0,T)$.

We also introduce spaces $\CC^k_{PC}$ consisting of $u\in \CC^k$
such that either $u\in \CC^k$ or there exists $\t=\t(u)\in [0,T]$
such that $\|u(\cdot,t)\|_{Z_T^k}$ is bounded,  $u(\cdot,t)$ is
continuous in $Z_T^k$ in $t\in[0,\t]$, and
 $u(\cdot,t)$ is continuous in $Z_T^k$ in $t\in[\t+\e,T]$ for any $\e>0$.
We also introduce spaces $Y^k_{PC}=X^{k}\!\cap \CC_{PC}^{k-1}$, with the norms from $Y^k$.

\subsection*{Conditions for the coefficients}
 To proceed further, we assume that Conditions
\ref{cond3.1.A}-\ref{condK} remain in force throughout this paper.
 \begin{condition} \label{cond3.1.A} The matrix  $b=b^\top$ is
symmetric  and bounded. In addition, there exists a constant
$\d>0$ such that
\be
 \label{Main1} y^\top  b
(x,t,\o)\,y-\frac{1}{2}\sum_{i=1}^N |y^\top\b_i(x,t,\o)|^2 \ge
\d|y|^2 \quad\forall\, y\in \R^n,\ (x,t)\in  D\times [0,T],\
\o\in\O. \ee
\end{condition}
\begin{condition}\label{cond3.1.B}
The functions  $f(x,t,\o)$, $\lambda (x,t,\o)$, $\b_i(x,t,\o)$,
and $\oo\b_i(x,t,\o)$, are bounded.
\end{condition}
\begin{condition}\label{condK} The mapping $\G: Y^1_{PC}\to Z_T^0$ is linear and
continuous.
\end{condition}
\par
Condition \ref{condK} allows, for instance, to consider $\G$ such
$\G u=u(\cdot,0)$, i.e., it covers periodic boundary value
conditions (\ref{period}).  Another example includes the case where
 there exists an integer $m\ge 0$, a set
$\{t_i\}_{i=1}^m\subset[0,T)$, and linear continuous operators
 $\ww\G_0:L_2([0,T];{\cal B}_1,\ell_1,H^0)\to
H^0$, $\ww\G_i:H^0\to H^0$, $i=1,..,N$, such that
\baaa \G u=\ww\G_0 u+\sum_{i=1}^m \ww\G_iu(\cdot,t_i). \label{G0Gi}\eaaa
In particular, it includes
 $$
\ww\G_0u=\int_0^Tk_0(t)u(\cdot,t)dt,\quad
\ww\G_iu(\cdot,t_i)=k_iu(\cdot,t_i),$$ where $k_0(\cdot)\in
L_2(0,T)$ and $k_i\in\R$. It covers also $\G$ such that $$
\ww\G_0u=\int_0^Tdt\int_Dk_0(x,y,t)u(y,t)dx,\quad
\ww\G_iu(\cdot,t_i)(x)=\int_Dk_i(x,y)u(y,t_i)dy,$$
 where $k_i(\cdot)$ are some
regular enough kernels.
\par We introduce the set of parameters $$
\ba {\cal P} \defi \biggl( n,\,\, D,\,\, T,\,\, \G,\
\delta,\,\,\,\,\\ \esssup_{x,t,\o,i}\Bigl[| b(x,t,\o)|+
|{f(x,t,\o)}|+|\lambda(x,t,\o)|+|\b_i(x,t,\o)|+|\oo\b_i(x,t,\o)|\Bigr].
\ea $$
\par
Sometimes we shall omit $\o$.
\subsection*{The definition of solution}
\begin{proposition} %2.2
\label{propL} Let $\zeta\in X^0$,
 let a sequence  $\{\zeta_k\}_{k=1}^{+\infty}\subset
L^{\infty}([0,T]\times\O, \ell_1\times\P;\,C(D))$ be such that all
$\zeta_k(\cdot,t,\o)$ are progressively measurable with respect to
$\F_t$, and let $\|\zeta-\zeta_k\|_{X^0}\to 0$. Let $t\in [0,T]$ and
$j\in\{1,\ldots, N\}$ be given.
 Then the sequence of the
integrals $\int_0^t\zeta_k(x,s,\o)\,dw_j(s)$ converges in $Z_t^0$ as
$k\to\infty$, and its limit depends on $\zeta$, but does not depend
on $\{\zeta_k\}$.
\end{proposition}
\par
{\it Proof} follows from completeness of  $X^0$ and from the
equality
\begin{eqnarray*}
\E\int_0^t\|\zeta_{k}(\cdot,s,\o)-\zeta_m(\cdot,s,\o)\|_{H^0}^2\,ds
=\int_D\,dx\,\E\left(\int_0^t\big(\zeta_k(x,s,\o)-
\zeta_m(x,s,\o)\big)\,dw_j(s)\right)^2.
\end{eqnarray*}
\begin{definition} %{2.1}
\rm Let $\zeta\in X^0$, $t\in [0,T]$, $j\in\{1,\ldots, N\}$, then we
define $\int_0^t\zeta(x,s,\o)\,dw_j(s)$ as the limit  in $Z_t^0$ as
$k\to\infty$ of a sequence $\int_0^t\zeta_k(x,s,\o)\,dw_j(s)$, where
the sequence $\{\zeta_k\}$ is such  as in Proposition \ref{propL}.
\end{definition}
\begin{definition} %3.1
\label{defsolltion} \rm Let $u\in Y^1$, $\chi_i\in X^0$,
$i=1,...,N$, and $\varphi\in X^{-1}$. We say that equations
(\ref{parab1})-(\ref{parab10}) are satisfied if \baaa
&&u(\cdot,t,\o)=u(\cdot,T,\o)+ \int_t^T\big(\A u(\cdot,s,\o)+
\varphi(\cdot,s,\o)\big)\,ds \ \nonumber
\\&&\hphantom{xxx}+ \sum_{i=1}^N
\int_t^TB_i\chi_i(\cdot,s,\o)ds-\sum_{i=1}^N
\int_t^T\chi_i(\cdot,s)\,dw_i(s)
%\eqno(3.2)
\label{intur} \eaaa for all $r,t$ such that $0\le r<t\le T$, and
this equality is satisfied as an equality in $Z_T^{-1}$.
\end{definition}
Note that the condition on $\p D$ is satisfied in the  sense that
$u(\cdot,t,\o)\in H^1$ for a.e. \ $t,\o$. Further, $u\in Y^1$, and
the value of  $u(\cdot,t,\o)$ is uniquely defined in $Z_T^0$ given
$t$, by the definitions of the corresponding spaces. The integrals
with $dw_i$ in (\ref{intur}) are defined as elements of $Z_T^0$. The
integral with $ds$ in (\ref{intur}) is defined as an element of
$Z_T^{-1}$. In fact, Definition \ref{defsolltion} requires for
(\ref{parab1}) that this integral must be equal  to an element of
$Z_T^{0}$ in the sense of equality in $Z_T^{-1}$.
\section{The main results}%{Problems with non-local conditions}
\begin{theorem}
\label{Th1}  There exist a number $\kappa=\kappa({\cal P})>0$ such
that problem (\ref{parab1})-(\ref{parab2}) has an unique solution
$(u,\chi_1,...,\chi_N)$ in the class $Y^1\times (X^0)^N$,  for any
$\varphi\in X^{-1}$, $\xi\in Z_T^0$, and any $\G$ such that $
\|\G\|\le \kappa$, where $\|\G\|$ is the norms of the operator $\G:
Y^1\to Z_T^0$. In addition, \be \label{3.3} \| u
\|_{Y^1}+\sum_{i=1}^N\|\chi_i\|_{X^0}\le C \left(\| \varphi \|
_{X^{-1}}+\|\xi\|_{Z_T^0}\right), \ee where $ C=C(\kappa,{\cal
P})>0$ is a constant that depends only on $\kappa$ and ${\cal P}$.
\end{theorem}
\par Let $\Ind$ denote  the indicator function.
\begin{theorem}
\label{Th2} Let $\G$ be such that
there exists $\t<T$ such that $\G u=\G (\Ind_{\{t\le
\t\}}u)$. Then  \be \label{3.4} \|
u\|_{Y^1}+\sum_{i=1}^N\|\chi_i\|_{X^0}\le C \left(\| \varphi \|
_{X^{-1}}+\|u\| _{X^{-1}}+ \|\xi\|_{Z_T^0}\right) \ee for all
solutions $(u,\chi_1,...,\chi_N)$ of problem
(\ref{parab1})-(\ref{parab2}) in the class $Y^1\times (X^0)^N$,
where $C=C({\cal P})>0$  depends only on ${\cal P}$ and does not
depend on $u$, $\varphi$ and  $\xi$.
\end{theorem}
\par
Starting from now and up to the end of this section, we assume that
Condition \ref{condB} holds.
\begin{condition}\label{condB}
\begin{enumerate}
\item
The domain $D$ is bounded. \item The functions  $b(x,t,\o)$,
$f(x,t,\o)$, $\lambda (x,t,\o)$, $\b_i(x,t,\o)$ and
$\oo\b_i(x,t,\o)$ are differentiable in $x$ for a.e. $t,\o$, and the
corresponding derivatives are bounded. \item $\b_i(x,t,\o)=0$ for
$x\in \p D$, $i=1,...,N$. \item $\F_0$ is the $\P$-augmentation of
the set $\{\emptyset,\O\}$.
\item
$\G=\G_0+\G_1$, where $\G_0:X^1\to Z_0^1$ and $\G_1:Y_{PC}^1\to
Z_0^1$ are continuous linear operators such that there exists $\t<T$
such that $\|\G_1(\Ind_{\{t\le\t\}}u)\|_{Z_0^1}\le \|u|_{t\le
\t}\|_{\CC^1(0,\t)}$ for all $u\in Y^1$  such that $u|_{t\le \t}\in
\CC^1(0,\t)$.
\end{enumerate}
\end{condition}
\par In particular, Condition \ref{condB}(ii) implies that there exist
modifications of $\b_i$ such that the functions $\b_i(x,t,\o)$ are
continuous in $x$ for a.e. $t,\o$. We assume that $\b_i$ are such
functions. \begin{example}{\rm The assumptions on $\G$ in Condition
\ref{condB} are satisfied, for instance, if there exists an integer
$m\ge 0$, a set $\{t_i\}_{i=1}^m\subset[0,T)$, and linear continuous
operators $\oo\G:L_2(Q)\to H^0$, $\oo\G_i:H^0\to H^0$, $i=0,1,..,N$,
such that the operators $\oo\G:L_2([0,T];{\cal B}_1,\ell_1,H^1)\to
W_2^1(D)$ and $\oo\G_i:H^1\to W_2^1(D)$ are continuous and
$$ \G u=\oo\G_0u(\cdot,0)+\E\{\oo\G
u+\sum_{i=1}^m \oo\G_iu(\cdot,t_i)\}. $$}
\end{example}
\begin{theorem}
\label{Th3}    Assume that problem (\ref{parab1})-(\ref{parab2})
with $\varphi\equiv 0$, $\xi\equiv 0$, does not admit non-zero
solutions
 $(u,\chi_1,...,\chi_N)$  in the class $Y^1\times
(X^0)^N$. Then   problem (\ref{parab1})-(\ref{parab2}) has a unique
solution $(u,\chi_1,...,\chi_N)$  in the class $Y^1\times (X^0)^N$,
 for any $\varphi\in X^{-1}$, and $\xi\in H^0$. In addition, \be
\label{3.5} \| u \|_{Y^1}+\sum_{i=1}^N\|\chi_i\|_{X^0}\le C \left(\|
\varphi \| _{X^{-1}}+\|\xi\|_{H^0}\right), \ee where $C>0$
 does not depend on $\varphi$ and  $\xi$.
\end{theorem}
\begin{theorem} \label{Th4} There exists $\e_0>0$ such that, for any $\e\in(-\e_0,\e_0)$
such that $\e\neq 0$, problem (\ref{parab1})-(\ref{parab10}) with
the boundary value condition \baa u(\cdot,T)-(1+\e)\G u=\xi
\label{ebv}\eaa has a unique solution $(u,\chi_1,...,\chi_N)$  in
the class $Y^1\times (X^0)^N$ for any $\varphi\in X^{-1}$ and
$\xi\in H^0$. In addition, \be \label{3.5e} \| u
\|_{Y^1}+\sum_{i=1}^N\|\chi_i\|_{X^0}\le C \left(\| \varphi \|
_{X^{-1}}+\|\xi\|_{H^0} \right), \ee where $C=C(\e)>0$
 does not depend on $\varphi$ and  $\xi$.
\end{theorem}
\begin{corollary}\label{corrAP} For any $\kappa\in R$,
there exists $\e_0=\e_0({\cal P},\kappa)>0$ such that, for any
$\e\in(-\e_0,\e_0)\backslash \{0\}$, for any $\varphi\in
X^{-1}$, there exists a unique solution $(u,\chi_1,...,\chi_N)$ in the class $Y^1\times
(X^0)^N$, of
problem (\ref{parab1}),(\ref{parab10}) with the boundary conditions \baa u(\cdot,T)=\k(1+\e)u(\cdot,0).
\label{ebp}\eaa
\end{corollary}
\begin{corollary}\label{corrAP2}
Corollary \ref{corrAP} with $\kappa=1$ implies that backward SPDE
(\ref{parab1})-(\ref{parab10}) can be regarded as an almost surely almost periodic solution on
a single time period.
\end{corollary}
Note that Tudor (1992) and other authors  considered  almost periodic solution for SPDEs on a
infinite time horizon, i.e., for many periods. We consider a single period only.
\begin{remark} The "almost periodicity" in the mean-square sense of  Tudor (1992)
requires that, for any $\e>0$, $\|u(\cdot,0)-u(\cdot,T)\|_{Z_T^0}\le
\e$; we establish that the equality $u(\cdot,0)=(1+\e)u(\cdot,T)$
can be achieved.  It follows that  $\|u(\cdot,0)-u(\cdot,T)\|_{Z_T^0} \le \e \|u(\cdot,T)\|_{Z_T^0}\le \e C\|\varphi\|_{X^{-1}}$.
Therefore, for a single time period, the "almost
periodicity" property in Corollary \ref{corrAP2}) implies the "almost periodicity" in the mean-square sense from Tudor (1992).  Moreover, this condition is
stronger than the mean-square almost periodicity since it requires  that the shapes of $u(\cdot,T)$ and $u(\cdot,0)$ are the same up
to proportionality.
\end{remark}
\par The remaining part of the section devoted to
an example where
an exact periodic condition is satisfied almost surely; this corresponds to the case  where
$\e=0$ in (\ref{ebp}).
\par
 Let functions ${\ww\b_i: Q\times
\O \to \R^n}$, $i=1,\ldots, M$, be such that $$
2b(x,t,\o)=\sum_{i=1}^N\b_i(x,t,\o)\,\b_i(x,t,\o)^\top
+\sum_{j=1}^M\,\ww\b_j(x,t,\o)\,\ww\b_j(x,t,\o)^\top, $$ and $\ww
\b_i$ has the similar properties as $\b_i$. (Note that, by Condition
\ref{cond3.1.A}, $2b>\sum_{i=1}^N\b_i\b_i^\top$).
\par
 Let
$\ww w(t)=(\ww w_1(t),\ldots, \ww w_M(t))$ be a new Wiener process
independent on $w(t)$. Let $a\in L_2(\O,\F,\P;\R^n)$ be a vector
such that $a\in D$. We assume also that $a$ is independent from
$(w(t)-w(t_1),\w w(t)-\w w(t_1))$ for all $t>t_1>s$. Let $s\in[0,T)$
be given. Consider the following Ito equation
\begin{eqnarray}
%\begin{array}{c}
\label{yxs} &&dy(t) = \ww
f(y(t),t)\,dt+\sum_{i=1}^N\b_i(y(t),t)\,dw_i(t) +\sum_{j=1}^M\ww
\b_j(y(t),t)\,d \ww w_j(t),
%\vspace*{20pt}
\nonumber\\ [-6pt] &&y(s)=x.
%\eqno{3.12}
%\end{array}
\end{eqnarray}
Here, $\ww f=\w f-\sum_i\oo\b_i\b_i$, $\ww f:D\times [0,T]\times \O
\to \R^n$ is a vector functions with the components $\ww f_i$.
\par
Let  $y(t)=y^{a,s}(t)$ be the solution of (\ref{yxs}), and let
$\tau^{a,s}\defi\inf\{t\ge s:\ y^{a,s}(t)\notin D\}$.

The following lemma is a modification  for the case of random coefficients of Lemma 2.1 from Dokuchaev (2004).
\begin{lemma}\label{propnu} There exists
$\nu\in(0,1)$ that depends only on  ${\cal P}$ such that
$\P(\tau^{a,0}> T )\le \nu$ for any random vector $a$ such that
$a\in D$ a.s. and $a$ does not depend on $w(t)-w(r)$ for all
$t>r>0$.
\end{lemma}
\begin{theorem}
\label{Th5} Let the functions $b,f$ and $\lambda$ be such that the
operator $\A$ can be represented as $$ \A
v=\sum_{i,j=1}^nb_{ij}(x,t,\o)\frac{\p^2 v}{\p x_i \p x_j}(x)
+\sum_{i=1}^n\w f_i(x,t,\o)\frac{\p v}{\p x_i
}(x)-\w\lambda(x,t,\o)v(x),
$$ where $\w\lambda(x,t,\o)\ge 0$ a.e., and where $\w f_i$ are bounded
functions.
\par
Further, let at least one of the following conditions is satisfied:
\begin{itemize}
\item[(i)] There exists  $c_\lambda>0$ such that
$\w\lambda(x,t,\o)\ge c_\lambda$ for all $x,t,\o$; or
\item[(ii)] $\kappa\in(-1,1)$; or
\item[(iii)]  For
$\nu\in (0,1)$  from Lemma \ref{propnu}, \baa
\frac{1}{2}\sum_{i=1}^N\int_0^T
\sup_{x,\o}\oo\b_i(x,t,\o)^2dt+\log\nu < 0. \label{smallb}\eaa
\end{itemize}
Furthermore, let $b\in \X_c^3$, $\w f\in\X_c^2$, $\w\lambda\in\X^1_c$, $
\b_i\in\X_c^3$. Then there exists $\oo\kappa>1$ such that, for any $\kappa\in[-\oo\kappa,\oo\kappa]$,
 problem
(\ref{parab1})-(\ref{parab10})  with the boundary condition
\baa u(\cdot,T)-\kappa u(\cdot,0)=\xi \label{kappa} \eaa  has
 a unique solution
$(u,\chi_1,...,\chi_N)$ in the class $Y^1\times(X^0)^N$  for any
$\varphi\in X^{-1}$ and $\xi\in Z^0_T$.
 In addition, (\ref{3.5})
holds with $C>0$ that
 does not depend on $\varphi$ and  $\xi$.
\end{theorem}
\index{For deterministic parabolic equations, similar conditions were
introduced in Dokuchaev (1994) (see also Theorem 2.2 from Dokuchaev
(2004)).}
\section{Proofs}
%\subsection*{Problems for forward equations}
Let $s\in (0,T]$, $\varphi\in X^{-1}$ and $\Phi\in Z^0_s$. Consider
the problem \be \label{4.1}
\begin{array}{ll}
d_tu+\left( \A u+ \varphi\right)dt +
\sum_{i=1}^NB_i\chi_i(t)dt=\sum_{i=1}^N\chi_i(t)dw_i(t), \quad t\le s,\\
u(x,t,\o)|_{x\in \p D}, \\
 u(x,s,\o)=\Phi(x,\o).
\end{array}
 \ee
 \par
The following lemma represents an analog of the so-called "the first
energy inequality", or "the first fundamental inequality" known for
deterministic parabolic equations (see, e.g., inequality (3.14) from
Ladyzhenskaya (1985), Chapter III).
\begin{lemma}
\label{lemma1} Assume that Conditions \ref{cond3.1.A}--\ref{condK}
are satisfied.  Then problem (\ref{4.1}) has an unique solution a
unique solution $(u,\chi_1,...,\chi_N)$ in the class
$Y^1\times(X^0)^N$  for any $\varphi\in X^{-1}(0,s)$, $\Phi\in
Z_s^0$, and \be \label{4.2} \| u
\|_{Y^1(0,s)}+\sum_{i=1}^N\|\chi_i\|_{X^0}\le C \left(\| \varphi \|
_{X^{-1}(0,s)}+\|\Phi\|_{Z^0_s}\right), \ee where $C=C({\cal P})$
does not depend on $\varphi$ and $\xi$.
\end{lemma}
(See, e.g., Dokuchaev (1991) or Theorem 4.2 from Dokuchaev (2010)).
\par Note that the solution $u=u(\cdot,t)$
is continuous in $t$ in $L_2(\O,\F,\P,H^0)$, since
$Y^1(0,s)=X^{1}(0,s)\!\cap \CC^{0}(0,s)$.
\par
Introduce  operators $L_s:X^{-1}(0,s)\to Y^1(0,s)$ and
$\L_s:Z^0_s\to Y^1(0,s)$, such that $u=L_s\varphi+\L_s\Phi,$ where
 $(u,\chi_1,...,\chi_N)$ is the solution of  problem (\ref{4.1})  in the class
$Y^2\times(X^1)^N$. By Lemma \ref{lemma1}, these linear operators
are continuous.
\par
Introduce   operators $\Q:Z_T^0\to Z_T^0$ and $\TT:X^{-1}\to Z_T^0$
such that  $\Q \Phi=\G\L_T\Phi$ and $\TT\varphi=\G L_T\varphi$,
i.e., $\Q\Phi+\TT\varphi= \G u$,
    where
$u$ is the solution in $Y^1$ of   problem (\ref{4.1}) with $s=T$,
$\varphi\in X^{-1}$, and $\Phi\in Z_T^0$.
It is easy to see that if the operator $\G:Y^1\to Z_T^0$ is continuous, then
the operators   $\Q:Z_T^0\to Z_T^0$ and $\TT:X^{-1}\to Z_T^0$ are linear and continuous. In particular,
$\|\Q\|\le \|\G\|\|\L_T\|$, where $\|\Q\|$, $\|\G\|$, and
$\|\L_T\|$, are the norms of the operators $\Q: Z_T^0\to Z_T^0$,
$\G: Y^1\to Z_T^0$, and $\L_T: Z_T^0\to Y^1$, respectively.

\begin{lemma}\label{lemmaQ} Assume that the operator $\G:Y^1\to Z_T^0$ is continuous.
 If the operator $(I-\Q)^{-1}:Z_T^0\to Z_T^0$ is also continuous
then problem (\ref{4.1})   has a unique solution
$(u,\chi_1,...,\chi_N)$ in the class $Y^1\times(X^0)^N$ for any
$\varphi \in X^{-1}$, $\Phi \in Z_T^0$. For this solution, \baa
\label{Q} u=L_T\varphi+\L_T (I-\Q)^{-1}(\xi+\TT\varphi) \eaa
 and \baaa
\label{Qes} \| u \|_{Y^1(0,s)}+\sum_{i=1}^N\|\chi_i\|_{X^0}\le C
\left(\| \varphi \| _{X^{-1}(0,s)}+\|\Phi\|_{Z^0_s}\right), \eaaa
where $C=C({\cal P})$ does not depend on $\varphi$ and $\xi$.
\end{lemma}
\par
{\it Proof of Lemma \ref{lemmaQ}}. For brevity, we denote
$u(\cdot,t)=u(x,t,\o)$. Clearly,
 $u\in Y^1$ is the solution of   problem
(\ref{parab1})-(\ref{parab2}) with some $(\chi_1,...,\chi_N)\in
(X^0)^N$ if and only if \baa
&&u=\L_Tu(\cdot,T)+L_T\varphi, \label{Q1}\\
&&u(\cdot,T)-\G u=\xi. \label{Q2}\eaa
 Since $ \G u=\Q
u(\cdot,T)+\TT\varphi$, equation (\ref{Q2})  can be rewritten as
\baa
 u(\cdot,T)-\Q u(\cdot,T)-\TT\varphi
=\xi. \label{Q3}\eaa  By the continuity of $(I-\Q)^{-1}$, equation
(\ref{Q3}) can be rewritten as
$$ u(\cdot,T)=(I-\Q)^{-1}(\xi +\TT\varphi). $$  Therefore, equations (\ref{Q1})-(\ref{Q2}) imply that \baaa \label{4.3}
u=L_T\varphi+\L_T u(\cdot,T)=L_T\varphi+\L_T
(I-\Q)^{-1}(\xi+\TT\varphi). \eaaa Further, let us show that if
(\ref{Q}) holds then equations (\ref{Q1})-(\ref{Q2})  hold. Let $u$
be defined by (\ref{Q}). Since $u=L_T\varphi+\L_T u(\cdot,T)$, it
follows that $u(\cdot,T)=(I-\Q)^{-1}(\xi+\TT\varphi)$. Hence \baaa
u(\cdot,T)-\Q u(\cdot,T)=\xi+\TT\varphi,\eaaa i.e., $u(\cdot,T)-\G
\L_T u(\cdot,T)=\xi+\TT\varphi=\xi+\G L_T\varphi.$   Hence \baaa
u(\cdot,T)-\G [\L_T u(\cdot,T)+L_T\varphi]=\xi.\eaaa This means that
(\ref{Q1})-(\ref{Q2}) hold. Then the proof of Lemma \ref{lemmaQ}
follows. $\Box$
 \par
{\it Proof of Theorem \ref{Th1}}. Since the operator $\Q: Z_T^0\to
Z_T^0$ is continuous, the operator $(I-\Q)^{-1}:Z_T^0\to Z_T^0$ is
continuous for small enough $\|\Q\|$, i.e. for a small enough
$\kappa>0$. Then the proof of Theorem \ref{Th1} follows. $\Box$
\par
{\it Proof of Theorem  \ref{Th2}}. For a real $q>0$, set
$u_q(x,t,\o)\defi e^{\q(T-t)}u(x,t,\o)$. Then  $u_q$ is the solution
of
  problem (\ref{parab1})-(\ref{parab2}) with  $\varphi$  replaced  by
$e^{\q(T-t)}\varphi(x,t)+q u_q(x,t)$, and with  $\G$ replaced by the
operator defined such that $\G_q$, where$$ \qquad \G_q u_q=\G u.$$
By the assumptions on $\G$, we have that \baaa \|\G_qu_q\|_{Z_T^0}=
\|\G u\|_{Z_T^0}= \|\G\Ind_{t\le\t}
u\|_{Z_T^0}=\left\|\G\Ind_{\{t\le\t\}}
e^{-\q(T-t)}u_q\right\|_{Z_T^0} \le e^{-q(T-\t)}\|\G\|\|
u\|_{X^0(0,\t)}.\eaaa It follows that $\|\G_q\|\to 0$ as $q\to
+\infty$, for the norm of the operator $\G_q:Y^1\to Z_T^T$. By Lemma
\ref{lemma1} and Theorem \ref{Th1}, it follows that, for a large
enough $q>0$, \baaa \| u_q \|_{Y^1}+\sum_{i=1}^N\|\chi_i\|_{X^0} \le
C_1 \left(\| e^{\q(T-t)}\varphi
+qu_q\| _{X^{-1}} +\|\xi\|_{H_0}\right)\\
\le C_2 \left(\| \varphi \| _{X^{-1}}+\|u_q\| _{X^{-1}}
+\|\xi\|_{H_0}\right), \eaaa where $C_1=C_1({\cal P})>0$ and
$C_2=C_2(q,{\cal P})>0$ do not depend on $u,\varphi,\xi$. Then the
proof of Theorem \ref{Th2} follows. $\Box$
\par
Starting from now, we assume that Condition \ref{condB} is satisfied, in addition
to Conditions \ref{cond3.1.A}-\ref{condK}.
\par
The following lemma represents an analog of the so-called "the
second energy inequality", or "the second fundamental inequality"
known for the deterministic parabolic equations (see, e.g.,
inequality (4.56) from Ladyzhenskaya (1985), Chapter III).
\begin{lemma}\label{Th3.1.1}
Problem (\ref{4.1})   has a unique solution $(u,\chi_1,...,\chi_N)$
in the class $Y^2\times(X^1)^N$ for any $\varphi \in X^0$, $\Phi \in
Z_T^1$, and \be \| u \|_{{Y}^2}+\sum_{i=1}^N\|\chi_i\|_{X^1}
 \le   C\left(  \|\varphi \|_{X^0} +
\|\Phi \|_{Z_T^1} \right), \label{3.1.3} \ee where $C>0$ does not
depend on $\varphi$ and $\Phi$; it depends on ${\cal P}$ an on the
supremums of the derivatives listed in Condition \ref{condB}(ii).
\end{lemma}
\par
The lemma above represents a reformulation of Theorem 3.4 from
Dokuchaev (2010) or Theorem 4.3 from Dokuchaev (2012a). In the cited
paper, this result was obtained under some strengthened version of
Condition \ref{cond3.1.A}; this was restrictive. In Du and Tang
(2012), this result was obtained without this restriction, i.e.,
under Condition \ref{cond3.1.A} only.
\par
\begin{lemma}
\label{lemma3} The operator $\Q:Z_T^0\to Z_T^0$ is compact.
\end{lemma}
\par
{\it Proof of Lemma  \ref{lemma3}}. Let $u=\L_0\Phi$, where $\Phi\in
Z_T^0$. By the semi-group property of backward SPDEs from Theorem
6.1 from Dokuchaev (2010), we obtain that
$u|_{t\in[0,s]}=\L_su(\cdot,s)$ for all $s\in (0,T]$. By Lemmas
\ref{lemma1} and \ref{Th3.1.1}, we have that \baaa \|\G_1
u\|^2_{Z_0^1}\le C_0
\sup_{\tau\in[0,\t]}\|u(\cdot,\tau)\|^2_{Z^1_0}\le
C_1\sup_{\tau\in[0,T]}\inf_{t\in[\tau,T]} \|
u(\cdot,t)\|^2_{Z_t^1}\\\le \sup_{\tau\in[0,T]}\frac{C_1}{T-\tau}
\int_{\tau}^T\|u(\cdot,t)\|_{Z_t^1}^2dt\le
\frac{C_2}{T-\t}\|\Phi\|^2_{Z_T^0}\eaaa and \baaa
 &&\|\G_0u\|^2_{Z_0^1}\le %%%%ushlo s H^1
C_3\E\int_0^T\|u(\cdot,t)\|^2_{Z^1_T}dt\le C_4 \|\Phi\|_{Z_T^0}.
 \eaaa for constants $C_i>0$ which
do not  depend on $\Phi$. Hence the operator  $\Q:Z_T^0\to Z_T^0$ represents a linear continuous
operator $\Q:Z_T^0\to W_2^1(D)$. Note that, by the definitions, $Z_T^0=L_2(\O,\F_0,\P,L_2(D))$. Since $\F_0$ is a trivial $\s$-algebra, the
convergence in $Z_T^0$ is equivalent to convergence in $H^0=L_2(D))$. Since the embeddings of $W_2^1(D)$ into $H^0$ and into
$Z_T^0$ are compact operators, the proof of Lemma \ref{lemma3}
follows. $\Box$
\par
{\it Proof of  Theorem \ref{Th3}}. By the assumptions, the equation
$\Q \Phi=\Phi$ has the only solution $\Phi=0$ in $H^0$. By Lemma
\ref{lemma3} and by the Fredholm Theorem, the operator
$(I-\Q)^{-1}:H^0\to H^0$ is continuous. Then the proof of  Theorem
\ref{Th3} follows from representation (\ref{4.3}). $\Box$
\par
{\it Proof of  Theorem \ref{Th4}}.  By Lemma \ref{lemma3} and by the
Fredholm Theorem again, for any $\e_0\in(0,1)$, there exists a
finite set $\Lambda\subset\C$ such that the operator $(\lambda
I-\Q)^{-1}:H^0\to H^0$ is continuous for all $\lambda\in
(1-\e_0,1+\e_0)\backslash \Lambda$. Then the proof of Theorem
\ref{Th4} follows from representation (\ref{4.3}) again. $\Box$
\par
{\it Corollary \ref{corrAP}} is a special case of  Theorem \ref{Th4}
with $\G u=u(\cdot,0)$.
\par {\it Proof of Lemma
\ref{propnu}.}  This proof represents a modification of the proof
for Lemma 2.1 from Dokuchaev (2004) for the case of random
coefficients.  For the case of $\oo\b=0$, the proof of Lemma
\ref{propnu} can be found in Dokuchaev (2012c).

 Let $\mu =(\ww f,\b,x,s)$.\par Clearly, there exists
a finite interval $D_1\defi (d_1,d_2)\subset \R$ and a bounded
domain $D_{n-1}\subset \R^{n-1}$ such that $D\subset D_1\times
D_{n-1}$.
\par
For $(x,s)\in D\times [0,T)$, let $\tau_1^{x,s}\defi\inf\{t\ge s: \
y^{x,s}_1(t)\notin D_1\}$, where $y^{x,s}_1(t)$ is the first
component of the vector
$y^{x,s}(t)=(y^{x,s}_1(t),...,y^{x,s}_n(t))$. We have that
\be\label{dd1} \P(\tau^{x,0}>T )\le \P(\tau_1^{x,0}>T
)=\P(y_1^{x,0}(t)\in D_1\ \forall t\in[0,T]). \ee
\par
Let \baaa M^{\mu}(t)\defi
\sum_{k=1}^N\int_s^th_k(y^{x,0}(r),r)dw_i(r)+\sum_{k=N+1}^{N+M}\int_s^th_k(y^{x,0}(r),r)d\ww
w_i(r),\quad t\ge s, \eaaa where $h=(h_1,..,h_{N+M})$ is a vector
 that represents the first
row of the matrix \baaa (\b_1,...,\b_N,\w\b_1,...,\w\b_M) \eaaa with
the values in $\R^{n\times (N+M)}$.
\par
 Let $\w D_1\defi (d_1+K_1,d_2+K_2)$, where
$K_1\defi -d_2-\vartheta \sup_{x,t,\o}|\w f_1(x,t,\o)|$, $K_2\defi
-d_1+\vartheta \sup_{x,t}|\w f_1(x,t,\o)|$. Clearly, $\w D_1$
depends only on $n,D$, and $c_f$. It is easy to see that
\be\label{dd2} \P(y_1^{x,0}(t)\in D_1\ \forall t\in[0,T ])\le
\P(M^{\mu}(t)\in \w D_1\ \forall t\in[0,T ]). \ee Further,
\be\label{dd3}h(y^{x,0}(t),t)^\top h(y^{x,0}(t),t)=
|h(y^{x,0}(t),t)|^2\in [\d,c_{\b}],\ee where \baaa \d=\inf_{x,0,\o,\
\xi\in\R^n:\ |\xi|=1}2\xi^\top b(x,t,\o)\xi,\quad
c_\b=\sup_{x,0,\o,\ \xi\in\R^n:\ |\xi|=1}2\xi^\top
b(x,t,\o)\xi.\eaaa Clearly, $M^{\mu}(t)$ is a martingale vanishing
at $s$ with quadratic variation process
$$[M^{\mu}]_t\defi \int_0^t|h(y^{x,0}(r),r)|^2dr,\qquad t\ge 0.
$$
\par Let $\t^{\mu}(t)\defi \inf\{r\ge 0:\ [M^{\mu}]_r>t\}$.
 Note that $\t^{\mu}(0)=0$, and
the function $\t^{\mu}(t)$ is  strictly increasing in $t$ given $x$.
By Dambis--Dubins--Schwarz Theorem (see, e.g., Revuz and Yor
(1999)), the process $B^{\mu}(t)\defi M(\t^{\mu}(t))$ is a Brownian
motion vanishing at $t=0$, i.e.,   $B^{\mu}(0)=0$, and
$M^{\mu}(t)=B^{\mu}([M^{\mu}]_t)$. Clearly, \be\label{dd4} \ba
\P(M^{\mu}(t)\in \w D_1\ \hphantom{x}\forall t\in[0,s+T
])&=\P(B^{\mu}([M^{\mu}]_t)\in \w D_1\ \hphantom{x}\forall
t\in[0,T ])\\
&\le \P(B^{\mu}(r)\in \w D_1\ \hphantom{x}\forall
r\in[0,[M^{\mu}]_{T}]). \ea \ee By (\ref{dd3}), $[M^{\mu}]_{T }\ge
\d T $ a.s. for all $x$. Hence \be\label{dd5}\P(B^{\mu}(r)\in \w
D_1\ \hphantom{x}\forall r\in[0,[M^{\mu}]_{T }])\le\P(B^{\mu}(r)\in
\w D_1\ \hphantom{x}\forall r\in[0,\d T ]). \ee By
(\ref{dd1})--(\ref{dd2}) and (\ref{dd4})--(\ref{dd5}), it follows
that
$$
\ba \sup_{\mu}\P(\tau^{x,0}>T )\le \nu\defi \sup_{\mu}
\P(B^{\mu}(r)\in \w D_1\ \hphantom{x}\forall r\in[0,\d T ]), \ea
$$
and $\nu=\nu({\cal P})\in(0,1)$.
 This completes the
proof of Lemma \ref{propnu}.
 $\Box$
\par
\par
{\it Proof of  Theorem \ref{Th5}}. Let us introduce operators
   $$ \A^* v\defi\sum_{i,j=1}^n\frac{\p^2 }{\p x_i \p x_j}
\left(b_{ij}(x,t)v(x)\right)-\sum_{i=1}^n\frac{\p}{\p x_i }\left(
\ww f_i(x,t)v(x)\right)-\w\lambda(x,t)v(x)$$ and \baaa B_i^*v\defi
 -\sum_{k=1}^n \frac{\p }{\p x_k }\,\big(\beta_{ik}(x,t,\o)\,v(x))+
 \oo \beta_i(x,t,\o)\,v(x),\qquad i=1,\ldots ,N. \label{AB*}\eaaa
Here $b_{ij}$, $x_i$, $\b_{ik}$ are the components of $b$,  $\b_i$,
and $x$.
\par
Let $\rho\in Z_{s}^0$, and let $p=p(x,t,\o)$ be the solution of  the
problem \baaa &&d_tp=\A^* p\, dt +
\sum_{i=1}^NB^*_ip\,dw_i(t), \quad t\ge s,\nonumber\\
&&p|_{t=s}=\rho,\quad\quad p(x,t,\o)|_{x\in \p D}=0.\label{p}\eaaa
By Theorem 3.4.8 from Rozovskii (1990), this boundary value problem
has an unique solution $p\in Y^1(s,T)$.
Introduce an operator $\M_s: Z_s^0\to Y^1(s,T)$  such that
$p=\M_s\rho$, where $p\in Y^1(s,T)$ is the solution this boundary value problem..
\par
The following lemma  from Dokuchaev (2005) represents an analog of the so-called "the
second energy inequality", or "the second fundamental inequality"
known for the deterministic parabolic equations (see, e.g.,
inequality (4.56) from Ladyzhenskaya (1985), Chapter III).
\begin{lemma}\label{Th3.1f}
Problem (\ref{4.1})  has an unique solution $p\in {Y}^2$ for any
$\rho \in Z_s^1$, and \be \| p\|_{{Y^2(s,T)}}
 \le   C
\|\rho \|_{Z_s^1}, \label{3.1f} \ee where  $C>0$ does not depend on
$\rho$. This $C$ depends on ${\cal P}$ and on the supremums of the
derivatives in Condition \ref{condB}.
\end{lemma}
\par
By Theorem 4.2 from Dokuchaev (2010), we have that
$\kappa p(\cdot,T)=\Q^*\rho$, i.e., \baa
(\rho,\Q\Phi)_{Z_0^0}=(\rho,\kappa v(\cdot,0))_{Z_0^0}=
(p(\cdot,T),\kappa v(\cdot,T))_{Z_T^0}=
(\kappa p(\cdot,T),\Phi)_{Z_T^0}\label{dual}\eaa for $v=\L_T\Phi$. (See
also Lemma 6.1 from Dokuchaev (1991) and related results in Zhou
(1992)). \par Suppose that there exists $\Phi\in Z_T^0$ such that
$\kappa v(\cdot,0)=v(\cdot,T)$ for $v=\L_T\Phi$, i.e.,
$v(\cdot,0)=\Q\Phi=\Phi$.
Let us show that $\Phi=0$ in this case.
\par
\def\QQ{{\bf Q}}
Since $\Q\Phi\in Z_0^0$, it follows that $\Phi\in H^0=Z_0^0$. Let $p=\M_0\rho$ and
$\oo p(x,t,0)=\E p(x,t,\o)$ (meaning the projection from $Z_T^0$ on
$H^0=Z_0^0$). Introduce an operator $\QQ: H^0\to H^0$  such that
$\kappa\oo p(\cdot,T)=\QQ\rho$. By (\ref{dual}), the properties of $\Phi$ lead to the equality \baa
(\rho-\kappa p(\cdot,T),\Phi(\cdot,T))_{Z_T^0}=(\rho-\kappa \oo
p(\cdot,T),\Phi(\cdot,T))_{H_0}=0 \quad\forall \rho\in
H^0.\label{pPhi}\eaa It suffices to show that the set $\{\rho-\kappa \oo
p(\cdot,T)\}_{\rho\in H^0}$ is dense in $H^0$. For this, it suffices
to show that the equation $\rho-\QQ\rho=z$ is solvable in $H^0$ for
any  $z\in H^0$.
\par
Let us show that the operator $\QQ:H^0\to H^0$ is compact. Let $p$
be the solution of (\ref{p}). This means that $\kappa \E p(\cdot,T)=\QQ
\rho$. By Lemma \ref{Th3.1.1}, it follows that \be \label{2f}
\|p(\cdot,\tau)\|_{Z_\tau^1}\le C\|p(\cdot,s)\|_{Z_s^1},\quad
\tau\in[s,T], \ee where $C_*>0$ is a constant that does not depend
on $p$, $s$, and $\tau$.
\par
We have that $p|_{t\in[s,T]}=\M_s p(\cdot,s)$ for all $s\in
[0,T]$, and, for $\tau>0$, \baaa \|\oo p(\cdot,T)\|^2_{W_2^1(D)}&\le
& C_0\|p(\cdot,T)\|^2_{Z^1_T}\le C_1\inf_{t\in[0,T]}\|
p(\cdot,t)\|^2_{Z_t^1}\\ &\le& \frac{C_1}{T} \int_0^T\|
p(\cdot,t)\|^2_{Z_t^1}dt  \le \frac{C_2}{T}\|p\|^2_{X^1} \le
\frac{C_3}{T}\|\Phi\|_{H^0} \eaaa for constants $C_i>0$ that do not
depend on $\Phi$. Hence the operator $\QQ:H^0\to H^1$ is continuous.
The embedding of $H^1$ into $H^0$ is a compact operator (see, e.g.,
Theorem 7.3 from Ladyzhenskaia (1985), Chapter I).
\par
Let us show that if \baa\kappa \oo p(\cdot,T)=\kappa \E
p(\cdot,T)=\QQ\rho=p(\cdot,0)\label{pp}\eaa for some $\rho\in H^0$
then $\rho=0$.
\par
 Let $\rho\in H^0$ be such that $\rho\ge 0$ a.e. and $\int_{D}\rho(x)dx=1$.
 Let $a\in L_2(\O,\F,\P;\R^n)$ be independent from the process  $(w(\cdot),\w w(\cdot))$
  such that $a\in D$ a.s. and it
has the probability  density function $\rho$. Let $p=\M_0\rho$, and
 let $y^{a,0}(t)$ be
the solution of  Ito equation (\ref{yxs}) with the initial condition
$y(0)=a$.

For $t\ge s$, set   \baaa &&\g_M^{x,s}(t)
\defi\exp\biggl[\sum_{i=1}^N\int_s^t\oo \b_i(y^{x,s}(s),s)\,dw_i(s)
-\sum_{i=1}^N\frac{1}{2}\,\int_s^t\oo\b_i(y^{x,s}(s),s)^2\,ds\biggr],\quad
\\&& \g^{x,s}(t) \defi\exp\biggl[-\int_s^t \w
\lambda(y^{x,s}(t),t)\,dt
 \biggr] \g_M^{x,s}(t).
\eaaa By Theorem 6.1 from Dokuchaev (2011), for all bounded
functions $\Phi\in Z_T^0$ and $u=\L_T\Phi$, we have that
\baa\E\int_{D} p(x,T,\o)\Phi(x,\o)dx=\int_{D}
p(x,0)u(x,0)dx=\E\Ind_{\{ \tau^{a,0}\ge
T\}}\g^{a,0}(T)\Phi(y^{a,0}(T)) \quad\hbox{a.s.} \label{probd}\eaa
If $D=\R^n$ and $\Ind_{\{\tau^{a,0}\le T\}}\equiv 1$, then this
equality follows from Theorem 5.3.1 from Rozovskii (2001). Equality
(\ref{probd}) means that $p(x,T,\o)$ is the conditional (given
$\F_T$) probability density function of the vector $y^{a,0}(T)$ if
the process $y^{a,0}(t)$ is killed at $\p D$ and if it is
 killed inside $D$ with the rate of killing $\w\lambda$. In
 particular, it follows that $p(x,t,\o)\ge 0$ a.e. and
 \baaa\E\int_{D}
p(x,T,\o)dx= \E\Ind_{\{ \tau^{a,0}\ge T\}}\g^{a,0}(T).
\label{probd20}\eaaa
\par
 Assume first that $\w\lambda(x,t,\o)\ge c_\lambda> 0$ for all $x,t,\o$, i.e., that
condition (i) is satisfied. In this case, $0\le \g^{a,0}(T)\le \nu_1
\g_M^{a,0}(T)$,  where $\nu_1\defi e^{-c_\lambda T}$, $\nu_1\in (0,1)$. Hence
 \baa\E\int_{D}
p(x,T,\o)dx= \E\Ind_{\{ \tau^{a,0}\ge T\}}\g^{a,0}(T)\le \nu_1\E
\g_M^{a,0}(T)=\nu_1. \label{nu1}\eaa
\par
\index{Assume  now that a strengthened condition (ii) is satisfied.
Assume first that $\oo\b_i\equiv 0$ for all $i$, i.e., that a
version of the condition (ii) is satisfied.  By Lemma \ref{propnu},
it follows that \baaa \E\int_{D} p(x,T,\o)dx =\E\Ind_{\{
\tau^{a,0}\ge T\}}\g^{a,0}(T) \le \int_{D} \E p(x,T,\o)dx\le
\ww\nu<1. \label{probd3}\eaaa  Here $\ww\nu\in(0,1)$ is the same as
in Lemma \ref{propnu}.}

 Assume now that $|\kappa|<1$, i.e., that
condition (ii) is satisfied. This case can be reduced to the case of condition (i) as the following.
The problem  $u(x,t)$ can be replaced by the problem for $u_q(x,t)=u(x,t)e^{q(T-t)}$ with $q=T^{-1}\log|\kappa|<0$. The new boundary value condition for $u_q$ is $\kappa e^{-q T}u_q(\cdot,0)+u_q(\cdot,T)=0$, i.e.,
\baaa
\frac{\kappa}{|\kappa|}u_q(\cdot,0)+u_q(\cdot,T)=0.\eaaa  In the new equation
for  $u_q(x,t)$, the coefficient $\lambda$ has to be replaced by $\lambda-q$.
It follows that condition (i) with $c_\lambda=-q$ is satisfied for the new problem.

Further, assume that (\ref{smallb}) is satisfied, i.e. that
condition (iii) is satisfied. Let $p>1$ and $q>1$ be such that
$1/p+1/q=1$. By Lemma \ref{propnu}, we have that it that \baaa \E\int_{D}
p(x,T,\o)dx =\E\Ind_{\{ \tau^{a,0}\ge T\}}\g^{a,0}(T) \le \|\Ind_{\{
\tau^{a,0}\ge T\}} \|_{L_p(\O)}\|\g_M^{a,0}(T)\|_{L_q(\O)}\\ =\P(
\tau^{a,0}\ge T)^{1/p} \|\g_M^{a,0}(T)\|_{L_q(\O)}\le
\nu^{1/p} \|\g_M^{a,0}(T)\|_{L_q(\O)} .
\label{probd3}\eaaa  Clearly, we have that
\baaa \|\g_M^{a,0}(T)\|^q_{L_q(\O)}\le
\exp\left(\frac{q^2-q}{2}\sum_{i=1}^N\int_0^T
\sup_{x,\o}\oo\b_i(x,t,\o)^2dt\right). \eaaa Hence \baaa
\|\g_M^{a,0}(T)\|_{L_q(\O)}\le
\exp\left(\frac{q-1}{2}\sum_{i=1}^N\int_0^T
\sup_{x,\o}\oo\b_i(x,t,\o)^2dt\right). \eaaa Further, $(q-1)/q=1/p$
and $q-1=q/p$. Hence \baa \E\int_{D} p(x,T,\o)dx\le \nu_2(q),\label{nu2}\eaa
where \baaa \nu_2(q)\defi\nu^{1/p}\exp\left(\frac{1}{p}\left[\log\nu
+\frac{q}{2}\sum_{i=1}^N\int_0^T
\sup_{x,\o}\oo\b_i(x,t,\o)^2dt\right]\right). \eaaa By
(\ref{smallb}), there exists $q_0>1$ such that $\nu_2(q_0)<1$.
  \index{\begin{remark} It can be seen
from the proof that Theorem \ref{Th5} still holds if $c_\lambda=0$
and $\oo\b_i\neq 0$, if $\sup_{x,t,\o,i}|\oo\b_i(x,t,\o)|$ is small
enough such that
  $\nu^{1/p}\|\g_M^{a,0}(T)\|_{L_q(\O)}<1$  for some
$q>1$, where $\nu$ is the same as in Lemma \ref{propnu}, for
$1/p+1/q=1$. In other words, it is required that
$\|\g_M^{a,0}(T)\|_{L_q(\O)}$ is close enough to 1.
\end{remark}}
\par
Let
$\nu_*\defi\max(\nu_1,\nu_2(q_0))$. We have that $\nu_*\in (0,1)$. By (\ref{nu1}), (\ref{nu2}), and by the linearity of problem (\ref{p}), it follows that \baa\int_{D}
\E p(x,T,\o)dx\le \nu_*\int_{D} \rho(x)dx \label{probd4}\eaa for
all non-negative $\rho(x)$.
%\end{document}
\par
 Suppose  that
(\ref{pp}) holds for $\rho\in H^0$. Let \baaa \rho_+(x)\defi\max(0,
\rho(x)),\quad \rho_-(x)\defi\max(0,-\rho(x)).\eaaa Let $p_+$ and
$p_-$ be the solutions of (\ref{p}) with $s=0$ and with $\rho$
replaced by $\rho_{\pm}$ respectively. Let $\oo p_{\pm}(x,t)=\E
p_{\pm}(x,t,\o)$. By the definitions, \baaa \oo p_+(\cdot,T)\defi
\QQ \rho_+, \quad \oo p_-(\cdot,T)\defi \QQ \rho_-.\eaaa By
(\ref{probd}), it follows that $\oo p_{\pm}(x,T)\ge 0$ for a.e. $x$.
By (\ref{probd4}), it follows that \baa\int_{D} \oo
p_{\pm}(x,T)dx\le \nu_*\int_{D} \rho_{\pm}(x)dx. \label{ppnu} \eaa
 \par
 Let us select any $\oo\kappa>1$ such that $\oo\kappa\nu_*<1$.
 \par
 Let us assume first that
$\rho_+\neq 0$ and that $\kappa\in[0,1]$.
It follows that there exist a
measurable set $D_0\subset D$ such that  $\mes(D_0)>0$ and that
$\rho(x)>0$ and
 $\int_{D_0}\oo p_+(x,T)dx\le \nu_*\int_{D_0}\rho(x)dx$ for all $x\in D_0$. It follows that
$\int_{D_0}\oo p(x,T)dx=\int_{D_0}p_+(x,T)dx-\int_{D_0}p_-(x,t)dx\le \nu_*\int_{D_0}\rho(x)dx$. Therefore, $\kappa\oo  p(\cdot,T)\neq \rho(\cdot)$
in this case. Similarly,
we can show that $\kappa \oo p(\cdot,T)\neq \rho$ if $\rho_-\neq 0$ and $\kappa\in[0,1]$.

Further, let us assume that $\kappa\in[-\oo\kappa,0)$. Let $D_+=\{x: \rho(x)\ge 0\}$,  $D_-=\{x: \rho(x)<0\}$.
By the assumptions,
 \baaa\int_{D_+} \rho(x)dx=\kappa\int_{D_+}\oo p(x,T)dx>0, \qquad \int_{D_-} \rho(x)dx=\kappa\int_{D_-}
\oo p(x,T)dx<0. \label{pp1} \eaaa We have that \baaa && 0\le
-\int_{D_+}\oo p(x,T)dx\le -\nu_*\int_{D_-}\rho(x)dx,\qquad 0\le
\int_{D_-}\oo p(x,T)dx\le\nu_* \int_{D_+}\rho(x)dx. \eaaa
 Hence \baaa &&
-\int_{D_+}\oo p(x,T)dx\le -\nu_*\kappa\int_{D_-} \oo
p(x,T)dx,\qquad \int_{D_-}\oo p(x,T)dx\le\nu_*\kappa\int_{D_+} \oo
p(x,T)dx. \eaaa Hence \baa && \int_{D_+}\oo p(x,T)dx\ge
\nu_*\kappa\int_{D_-} \oo p(x,T)dx,\qquad \int_{D_-}\oo
p(x,T)dx\le\nu_*\kappa\int_{D_+} \oo p(x,T)dx. \label{pp2} \eaa
Hence \baaa \left|\int_{D_+}\oo p(x,T)dx\right|\le
|\nu_*\kappa|\left| \int_{D_-} \oo p(x,T)dx\right|,\qquad
\left|\int_{D_-}\oo p(x,T)dx\right|\le|\nu_*\kappa|\left|\int_{D_+}
\oo p(x,T)dx\right|. \label{pp3} \eaaa The system of the last two
inequalities can be satisfied only if the integrals there are zero.
This means that $\rho=0$.
\par
We have proved that  if (\ref{pp}) holds for $\rho\in H^0$ then
$\rho=0$. We had proved also that the operator $\QQ$ is compact. By
the Fredholm Theorem, it follows that the equation $\rho-\QQ\rho=z$
is solvable in $H^0$ for any $z\in H^0$. By (\ref{pPhi}), it follows
that $\Phi=0$.
 Therefore, the condition $\kappa u(\cdot,0)=u(\cdot,T)$
fails to be satisfied for $u\neq 0$, $\xi=0$, and $\varphi=0$. Thus,
$u=0$ is the unique solution of problem
(\ref{parab1})-(\ref{parab2}) for $\xi= 0$ and $\varphi=0$. Then the
proof of Theorem \ref{Th5} follows from Theorem \ref{Th3}.
$\Box$

\subsection*{Acknowledgment} This work  was
supported by ARC grant of Australia DP120100928 to the
author.
\section*{References} $\hphantom{XX}$
Al\'os, E., Le\'on, J.A., Nualart, D. (1999).
 Stochastic heat equation with random coefficients,
 {\it
Probability Theory and Related Fields} {\bf 115} (1), 41--94.

Arnold, L.,  and Tudor, C. (1998). Stationary and Almost Periodic
Solutions of Almost Periodic Affine Stochastic Differential
Equations, {\em Stochastics and Stochastic Reports} {\bf 64},
177--193.

\par
Bally, V., Gyongy, I., Pardoux, E. (1994). White noise driven
parabolic SPDEs with measurable drift. {\it Journal of Functional
Analysis} {\bf 120}, 484--510.

 Bedouhene, F., Mellah, O.,
Raynaud de Fitte. P. (2012). Bochner-Almost Periodicity for
Stochastic Processes. {\em Stochastic Analysis and Applications},
vol. 30, no. 2, pp. 322-342.

Bezandry, P. and Diagana, T. (2007). Square-mean almost periodic
solutions nonautonomous stochastic differential equations. {\em
Electron. J. Diff. Equ. } Vol. 2007, No. 117, pp. 1-10.

\par Caraballo,T., P.E. Kloeden, P.E.,
Schmalfuss, B. (2004). Exponentially stable stationary solutions for
stochastic evolution equations and their perturbation, {\em Appl.
Math. Optim.}, {\bf  50}, 183--207.
\par
 Chojnowska-Michalik, A. (1987). On processes of Ornstein-Uhlenbeck type in
Hilbert space. {\em Stochastics} {\bf  21}, 251--286.

\par
Chojnowska-Michalik, A. (1990). Periodic distributions for linear
equations with general additive noise, {\em  Bull. Pol. Acad. Sci.
Math.} {\bf 38} (1–12) 23--33.

\par
Chojnowska-Michalik, A., and Goldys, B. (1995). {Existence,
uniqueness and invariant measures for stochastic semilinear
equations in Hilbert spaces},  {\it Probability Theory and Related
Fields},  {\bf 102}, No. 3, 331--356.

Crewe, P.  (2013). Almost periodic solutions to stochastic evolution
equations on Banach spaces. {\em Stoch. Dyn.} {\bf 13}, 1250027, 23
pages.

Da Prato,G., and Tudor, C. (1995). Periodic and Almost Periodic
Solutions for Semilinear Stochastic Evolution Equations, {\em Stoch.
Anal. Appl.}  13 (1), 13–33.
\par
Da Prato, G., and Tubaro, L. (1996). { Fully nonlinear stochastic
partial differential equations}, {\it SIAM Journal on Mathematical
Analysis} {\bf 27}, No. 1, 40--55.

\par
Dokuchaev, N.G. (1992). { Boundary value problems for functionals
of
 Ito processes,} {\it Theory of Probability and its Applications}
 {\bf 36} (3), 459-476.
\index{\par
 Dokuchaev, N.G. (1994). Parabolic equations without the Cauchy
boundary condition and problems on control over diffusion processes.
I. {\em Differential Equations} {\bf 30}, No. 10, 1606-1617;
translation from Differ. Uravn. 30, No.10, 1738-1749.
\par
Dokuchaev, N.G. (1995). Parabolic equations without Cauchy boundary
condition and control problems
    for diffusion processes. II.
Differential  Equations {\bf 31}, No. 8, 1362-1372; translation from
Differ. Uravn. 31, No.8, 1409-1418.}
 \par Dokuchaev, N.G. (2004).
Estimates for distances between first exit times via parabolic
equations in unbounded cylinders. {\it Probability Theory and
Related Fields}, {\bf 129}, 290 - 314.
\par
Dokuchaev, N.G. (2005).  Parabolic Ito equations and second
fundamental inequality.  {\it Stochastics} {\bf 77} (2005), iss. 4.,
pp. 349-370.
\par
 Dokuchaev N. (2008). Parabolic Ito equations with mixed in time
conditions.
{\it Stochastic Analysis and Applications} {\bf 26}, Iss. 3, 562--576. %May 2008
\par
Dokuchaev, N. (2010). Duality and semi-group property for backward
parabolic Ito equations. {\em Random Operators and Stochastic
Equations. } {\bf 18}, 51-72.
\par
Dokuchaev, N. (2011). Representation of functionals  of Ito
processes in bounded domains. {\em Stochastics} {\bf 83}, No. 1,
45--66.
\par
 Dokuchaev, N. (2012a).
Backward parabolic Ito equations and second fundamental inequality.
{\em Random Operators and Stochastic Equations} {\bf 20}, iss. 1,
69-102.
\par
Dokuchaev, N. (2012b). Backward SPDEs with non-local in time and
space boundary conditions. Working paper, arXiv:1211.1460 (submitted).
\par
Dokuchaev, N. (2012c). On forward and backward SPDEs with
non-local boundary conditions.  Working paper in  arXiv (submitted).  
\par Du
K., and Tang, S. (2012). Strong solution of backward stochastic
partial differential equations in $C^2$ domains.  {\em Probability
Theory and Related Fields)} {\bf 154}, 255--285.
\par
Dorogovtsev, A.Ya., and  Ortega, O.A. (1988). On the Existence of
Periodic Solutions of a Stochastic Equation in a Hilbert Space.
{\em Visnik Kiiv. Univ. Ser. Mat. Mekh.} {\bf 115} No. 30,  21--30.

\par
Duan J.,  Lu K., Schmalfuss B. (2003). Invariant manifolds for
stochastic partial differential equations.{\em Ann. Probab.} {\bf
31}
 2109–2135.
\par
Feng C., Zhao H. (2012). Random periodic solutions of SPDEs via
integral equations and Wiener-Sobolev  compact embedding. {\em
Journal of Functional Analysis} {\bf 262}, 4377--4422.
\par
Fife, P.C. (1964). Solutions of parabolic boundary problems existing for all time. {\em Archive for Rational Mechanics and Analysis} 16
155–186.
\par
Gy\"ongy, I. (1998). Existence and uniqueness results for semilinear
stochastic partial differential equations. {\it Stochastic Processes
and their Applications} {\bf 73} (2), 271--299.
\par
Hess, P. (1991). Periodic Parabolic Boundary Value Problems and Positivity.     Wiley, New York.
\par
Kl\"unger, M. (2001). Periodicity and Sharkovsky's theorem for
random dynamical systems, {\em Stochastic and Dynamics} {\bf 1},
iss.3, 299--338.
\par Krylov, N. V. (1999). An
analytic approach to SPDEs. Stochastic partial differential
equations: six perspectives, 185--242, Mathematical Surveys and
Monographs, {\bf 64}, AMS., Providence, RI, pp.185-242.
\par
Ladyzhenskaia, O.A. (1985). {\it The Boundary Value Problems of
Mathematical Physics}. New York: Springer-Verlag.
\par
Lieberman G.M. (1999). Time-periodic solutions of linear parabolic differential equations. {\em Comm. Partial Differential
Equations} {\bf 24},  631–663.
\par
Liu, Y., Zhao, H.Z (2009). Representation of pathwise stationary
solutions of stochastic Burgers equations, {\em Stochactics and
Dynamics} {\bf  9} (4), 613--634.
\par
Maslowski, B. (1995). { Stability of semilinear equations with
boundary and pointwise noise}, {\it Annali della Scuola Normale
Superiore di Pisa - Classe di Scienze} (4), {\bf 22}, No. 1,
55--93.
\par
Mattingly. J. (1999). Ergodicity of 2D Navier–Stokes equations with random forcing and large viscosity. {\em Comm. Math.
Phys.} 206 (2),  273–288.
\par
 Mellah, O.,  and Raynaud de Fitte, P. (2007). Counterexamples to mean square almost periodicity of the solutions
 of some SDEs with almost periodic coefficients.
 {\em Electron. J. Differential Equations}, No. 117, 10 pp. (electronic).
 \par
Mohammed S.-E.A., Zhang T.,  Zhao H.Z. (2008). The stable manifold theorem for semilinear stochastic evolution equations
and stochastic partial differential equations. {\em Mem. Amer. Math. Soc.} 196 (917)  1–105.

\par Nakao, M. (1984).
Periodic solution of some nonlinear degenerate parabolic
equations. {\it Journal of Mathematical Analysis and
Applications} {\bf 104} (2), 554--557.
\par
Pardoux, E. (1993). Bulletin des Sciences Mathematiques, 2e Serie,
 {\bf 117}, 29-47.
\par
 Revuz, D., and Yor, M. (1999). {\it Continuous Martingales
and Brownian Motion}. Springer-Verlag: New York.
\par
Rodkina, A.E. (1992). On solutions of stochastic equations with
almost surely periodic
    trajectories.  {\em Differ. Uravn.} 28, No.3, 534--536 (in Russian).\par
Rozovskii, B.L. (1990). {\it Stochastic Evolution Systems; Linear
Theory and Applications to Non-Linear Filtering.} Kluwer Academic
Publishers. Dordrecht-Boston-London.
%315 p.
\par
Sinai, Ya. (1996). Burgers system driven by a periodic stochastic
flows, in: Ito's Stochastic Calculus and Probability Theory,
Springer, Tokyo, 1996, pp. 347–353.
\par
Shelukhin, V.V. (1993). Variation principle for non-local in time
problems for linear evolutionary equations. {\it Siberian Math. J.}
{\bf 34} (2), 191--207.
\par
Tudor, C. (1992).  Almost periodic solutions of affine stochastic
evolutions equations. {\em Stochastics and Stochastics Reports} {\em
38}, 251--266.
\par Vejvoda, O. (1982). Partial Differential Equations:
Time Periodic Solutions, Martinus Nijhoff, The Hague.
\par
Walsh, J.B. (1986). An introduction to stochastic partial
differential equations, Lecture Notes in Mathematics {\bf 1180},
Springer Verlag.
\par
Yong, J., and Zhou, X.Y. (1999). { Stochastic controls: Hamiltonian
systems and HJB equations}. New York: Springer-Verlag.
\par
 Zhou, X.Y. (1992). { A duality analysis on stochastic partial
differential equations}, {\it Journal of Functional Analysis} {\bf
103}, No. 2, 275--293.
\end{document}